\documentclass[final,a4paper]{siamltex}
\usepackage{amsmath,amssymb,latexsym}
\usepackage[final]{graphicx}
\newcommand{\even}{\mathop{\rm even}}
\newcommand{\odd}{\mathop{\rm odd}}
\newcommand{\field}[1]{\mathbb{#1}}
\newcommand{\R}{\field{R}}
\newcommand{\Z}{\field{Z}}
\newcommand{\N}{\field{N}}
\newcommand{\C}{\field{C}}
\renewcommand{\P}{\field{P}}

\renewcommand{\Re}{\mathop{\rm Re}}

\title{Orthogonality of Jacobi polynomials with general parameters}%

\author{A.B.J.\ Kuijlaars\thanks{Department of Mathematics,
Katholieke Universiteit Leuven (BELGIUM)}, corresponding author,
email: {\tt arno@wis.kuleuven.ac.be}
 \and A.\
Mart{\'\i}nez-Finkelshtein\thanks{University of Almer{\'\i}a and Instituto
Carlos I de F{\'\i}sica Te{\'o}rica y Computacional, Granada University
(SPAIN)}, email: {\tt andrei@ual.es}  \and R.\
Orive\thanks{University of La Laguna, Canary Islands (SPAIN)},
email: {\tt rorive@ull.es}
 }%
\date{\today}
\begin{document}
\maketitle

\begin{abstract}
In this paper we study the orthogonality conditions satisfied by
Jacobi polynomials $P_n^{(\alpha,\beta)}$ when the parameters
$\alpha$ and $\beta$ are not necessarily $>-1$. We establish
orthogonality on a generic closed contour on a Riemann surface. Depending
on the parameters, this leads to either full orthogonality
conditions on a single contour in the plane, or to multiple orthogonality
conditions on a number of contours in the plane.
In all cases we show that the orthogonality conditions characterize the
Jacobi polynomial $P_n^{(\alpha , \beta )}$ of degree $n$ up to a constant factor.
\end{abstract}
% ----------------------------------------------------------------



\pagestyle{myheadings} \thispagestyle{plain} \markboth{KUIJLAARS,
MART{\'I}NEZ-FINKELSHTEIN AND ORIVE}{ORTHOGONALITY FOR JACOBI
POLYNOMIALS}

\section{Introduction}

The Jacobi polynomials $P_n^{(\alpha,\beta)}$ are
given explicitly by
\begin{equation*}\label{explJac}
P_n^{(\alpha,\beta)} (z)=2^{-n} \sum_{k=0}^n \binom{n+\alpha}{n-k}
\binom{n+\beta}{k}(z-1)^k (z+1)^{n-k},
\end{equation*}
or, equivalently, by the well-known Rodrigues formula
\begin{equation}%\label{RodrJac}
\label{Rodrigues}
P_n^{(\alpha,\beta)} (z)=\frac{1}{2^n n!} (z-1)^{-\alpha}
(z+1)^{-\beta} \left( \frac{d}{dz} \right)^n \left[
(z-1)^{n+\alpha} (z+1)^{n+\beta}\right].
\end{equation}
These expressions show that $P_n^{(\alpha,\beta)}$ are  analytic
functions of the parameters $\alpha$ and $\beta$, and thus can
be considered for general $\alpha, \beta \in \C$.

The classical Jacobi polynomials correspond to parameters $\alpha, \beta>-1$.
For these parameters, the Jacobi polynomials are orthogonal on $[-1,1]$
with respect to the weight function $(1-z)^\alpha (1+z)^\beta$. As a result,
all their zeros are simple and belong to the interval $(-1,1)$. For
general $\alpha, \beta$, this is no longer valid. Indeed, the zeros
can be non-real, and there can be multiple zeros. In fact,
$P_{n}^{\left( \alpha ,\beta \right) }$ may have a multiple zero at
$z=1$ if $\alpha \in \{-1,\ldots,-n\}$, at $z=-1$ if $\beta \in
\{ -1,\ldots,-n\} $ or, even, at $z=\infty $ (which means a degree
reduction) if $n+\alpha +\beta \in \{-1, \ldots,-n\}$.

More precisely, for $k\in \{1,\ldots,n\}$, we have (see
\cite[formula (4.22.2)]{szego:1975}),
\begin{equation}
P_{n}^{\left( -k,\beta \right) }(z) = \frac{\Gamma(n+\beta+1)
}{\Gamma(n+\beta +1-k)}\, \frac{(n-k)!}{n!} \left(
\frac{z-1}{2}\right) ^{k}P_{n-k}^{\left( k,\beta \right) }\left(
z\right). \label{integer 1}
\end{equation}
This implies in particular that $P_{n}^{( -k,\beta) }(z) \equiv 0$ if additionally
$\max \left\{ k,-\beta \right\} \leq n\leq k-\beta -1$. Analogous
relations hold for $P_{n}^{\left( \alpha ,-l\right) }$ when $l\in
\{1,\ldots,n\}$. Thus, when both $k,l\in \mathbb{N}$ and $k+l\leq n
$, we have
\begin{equation}
P_{n}^{\left( -k,-l\right) }\left( z\right) =2^{-k-l}\left(
z-1\right) ^{k}\left( z+1\right) ^{l}P_{n-k-l}^{\left( k,l\right)
}\left( z\right)\,. \label{integer 2}
\end{equation}
Furthermore, when $n+\alpha +\beta  =-k \in \{-1, \ldots,-n\}$,
\begin{equation}
P_{n}^{\left( \alpha ,\beta  \right) }\left( z\right)
=\frac{\Gamma(n+\alpha+1) }{\Gamma(k+\alpha )} \,
\frac{(k-1)!}{n!}\, P_{k-1}^{\left( \alpha ,\beta \right) }\left(
z\right), \label{integer 3}
\end{equation}
see \cite[Eq.\ (4.22.3)]{szego:1975}; see \S 4.22 of \cite{szego:1975}
for a more detailed discussion. Formulas (\ref{integer
1})--(\ref{integer 3}) allow to exclude these special
integer parameters from our analysis.

In this paper we will show that for general $\alpha, \beta \in \C$,
but excluding some special cases, the Jacobi polynomials
$P_n^{(\alpha,\beta)}$ may still be characterized by orthogonality
relations. The case $\alpha, \beta >-1$ is classical, and standard
orthogonality on the interval $\left[ -1,1\right] $
takes place, being the key for the study of many properties of
Jacobi polynomials. Thus, our goal is to establish orthogonality
conditions for the remaining cases. We will show that
$P_n^{(\alpha, \beta)}$ satisfies orthogonality conditions on
certain curves in the complex plane. In some cases the orthogonality
conditions on a single curve are enough to characterize the Jacobi polynomial,
while in others a combination of orthogonality conditions on two or three
curves is required. This last phenomenon is called multiple orthogonality, see e.g.\
\cite{aptekarev:1998,nikishin/sorokin:1991}. In some particular cases this
orthogonality has been established before (see e.g.\
\cite{Arora:95}--\cite{Chen95}), and used in the study of
asymptotic behavior of these polynomials \cite{MR2002d:33017}.
Similar orthogonality conditions, but for Laguerre polynomials,
have been applied in \cite{MR1858305} in order to study the zero
distribution in the case of varying parameters, and in
\cite{Kuijlaars/Mclaughlin:01a}--\cite{Kuijlaars/Mclaughlin:01b},
in order to establish the strong asymptotics by means of the
Riemann-Hilbert techniques. A different kind of orthogonality
involving derivatives has been found for negative integer values
of the parameters of $P_n^{(\alpha , \beta )}$, see
\cite{alfaro/alvarez/rezola:2000}--\cite{alvarez-de-morales/etal:1998}.

We believe that these new orthogonality conditions can be useful
in the study of the zeros of Jacobi polynomials.
For general $\alpha, \beta \in \C$ the zeros are not confined to
the interval $[-1,1] $ but they distribute themselves in the complex
plane. K.\ Driver, P.\ Duren and collaborators
\cite{Driver/Duren99}--\cite{Duren/Guillou01} noted that the
behavior of these zeros is very well organized, see also
\cite{MR2002d:33017}. We believe that the orthogonality conditions
we find, and in particular the Riemann-Hilbert problem derived from
that (see Section 3 below) can be used to
establish asymptotic properties of Jacobi polynomials. In particular
this could explain the observed behavior of zeros.

\section{Orthogonality on a Riemann surface}

\begin{figure}[htb]
\centering \includegraphics[scale=0.6]{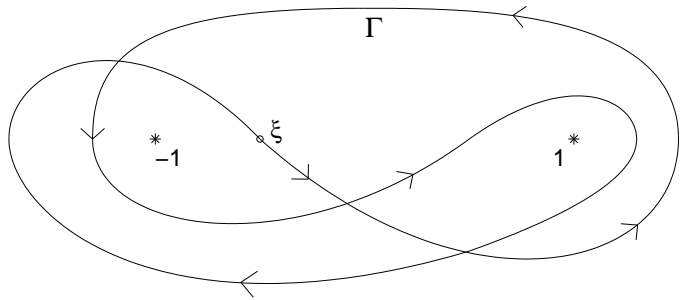} \caption{ Path
$\Gamma$.}\label{fig:pathGamma}
\end{figure}

Consider the path $\Gamma$ encircling the points $+1$ and $-1$
first in a positive sense and then in a negative sense, as shown
in Fig.~\ref{fig:pathGamma}. The point $\xi \in (-1,1)$ is the
begin and endpoint of $\Gamma$.

%%%%%%%%%%%%%%%%%%%%%%%%%%%%%%%%%%%%%%%%%%%%%%%%%%%%%%%%%%
%%%%%%%%%%%%%%%%%%%%%%%%%%%%%%%%%%%%%%%%%%%%%%%%%%%%%%%%%%

For $\alpha , \beta  \in \C$, denote
\begin{equation*}\label{weight}
w(z; \alpha , \beta ):=(1-z)^\alpha (z+1)^\beta=\exp [\alpha \log
(1-z) + \beta \log (z+1)].
\end{equation*}
It is a multi-valued function with branch points at $\infty$ and
$\pm 1$. However, if we start with a value of $w(z;\alpha, \beta)$ at a particular
point of $\Gamma$, and extend the definition of $w(z; \alpha, \beta)$
continuously along $\Gamma$, then we obtain a single-valued function
$w(z; \alpha, \beta)$ on $\Gamma$ if we view $\Gamma$ as a contour on
the Riemann surface for the function $w(z;\alpha, \beta)$.
For definiteness,
we assume that  the ``starting point'' is $\xi \in (-1,1)$, and
that the branch of $w$ is such that $w(\xi; \alpha , \beta )>0$.

In the sequel we prefer to view $\Gamma$ as a subset of the complex plane.
Then $\Gamma$ has points of self-intersection, see Fig.~\ref{fig:pathGamma}.
At points of self-intersection the value of $w(z;\alpha,\beta)$ is not
well-defined.

The following is our main result. It shows that the Jacobi polynomials
satisfy non-hermitian orthogonality conditions on $\Gamma$.
\begin{theorem} \label{thm:main}
Let $\alpha, \beta  \in \C $, and let $\Gamma$, $w(z; \alpha , \beta
)$ have the meaning as described above. Then for $k \in \{0, 1, \ldots, n\}$,
\begin{equation}\label{orthMain}
    \int_{\Gamma} t^k \, P_n^{(\alpha , \beta )}(t) w(t; \alpha,
    \beta)\,
    dt = \frac{-  \pi^2 2^{n+\alpha +\beta +3} e^{\pi i (\alpha +\beta
    )}}{\Gamma(2n+\alpha +\beta +2) \Gamma(-n-\alpha )
    \Gamma(-n-\beta)}\, \delta_{kn}\,.
\end{equation}
\end{theorem}
\begin{proof}
In the proof we use $f^{(k)}$ to denote the $k$-th derivative of $f$.

By the Rodrigues formula (\ref{Rodrigues}),
\begin{equation}\label{RodriguesBis}
P_n^{(\alpha,\beta)} (t)=\frac{(-1)^n}{2^n n!} \, \frac{w^{(n)}(t;
n + \alpha, n+\beta)}{ w(t; \alpha,  \beta)}.
\end{equation}
Integrating in (\ref{orthMain}) $n$ times by parts and using
(\ref{RodriguesBis}), we get
 \begin{equation}\label{byParts}
\begin{split}
  \int_{\Gamma} t^k P_n^{(\alpha , \beta )}(t)  w(t;\alpha, \beta)\,
 dt= \frac{(-1)^n}{2^n n!} \,
 \sum_{j=0}^{n-1} (-1)^j \big[t^k \big]^{(j)}    w^{(n-j-1)}(t; n+\alpha, n+\beta)
 \bigg|_{\Gamma} \\
 +\frac{1}{2^n n!} \, \int_{\Gamma} \big[ t^k \big]^{(n)} w(t; n+\alpha , n+\beta )\,
 dt\,.
\end{split}
\end{equation}
Since $w(z; \alpha, \beta)$ is single-valued on $\Gamma$,
\begin{equation*}\label{zeroMain}
\big[t^k \big]^{(j)}    w^{(n-j-1)}(t; n+\alpha, n+\beta)
\bigg|_{\Gamma}=0\,, \qquad \mbox{for } 0 \leq j \leq n-1\,.
\end{equation*}
Thus, if $k\leq n-1$, all the terms in the right-hand side of
(\ref{byParts}) vanish, which proves that the integral in
(\ref{orthMain}) is $0$ for $k=0, 1, \ldots , n-1$. Furthermore,
for $k=n$, we get
\[
 I_n(\alpha , \beta ):=\int_{\Gamma} t^n\, P_n^{(\alpha , \beta )}(t)  w(t;\alpha, \beta)\,
 dt=  2^{-n} \int_{\Gamma} w(t; n+\alpha, n+\beta)\, dt\,.
\]
Observe that $I_n(\alpha,\beta)$ is an analytic function of
$\alpha$ and $\beta$, so that we may compute it for a certain
range of parameters and then extend it analytically elsewhere.
Following \cite[\S 12.43]{Whittaker96}, we assume $\Re \alpha >
0$, $\Re \beta
>0$ and deform the path $\Gamma$, tautening it between $-1$ and
$+1$. Thus, $\Gamma$ will become the union of two small circles
around $\pm 1$ and
 two straight lines along $[-1,1]$, each piece traversed twice, once
in a positive direction and once in a negative direction.

Since the integrand is bounded in the neighborhoods of $\pm 1$,
$I_n(\alpha , \beta )$ splits in the following 4 integrals:
\[
I_n(\alpha , \beta )= \int_{-1}^1 f(t)\, dt - e^{2\pi
i(\alpha+n)}\, \int_{-1 }^{1} f(t)\, dt + e^{2\pi i(\alpha+\beta
+2 n)}\, \int_{-1 }^{1} f(t)\, dt- e^{2\pi i(\beta+n)}\, \int_{-1
}^{1} f(t)\, dt \,,
\]
where
$$
f(t)=2^{-n} w(t; n+\alpha, n+\beta)>0 \; \text{ for } \; t \in
(-1,1)\,.
$$
Thus,
\[
I_n(\alpha , \beta )= \left(1- e^{2\pi i(\alpha+n)} \right)\,
\left(1- e^{2\pi i(\beta+n)}\right)\, \int_{-1}^1 f(t)\, dt= -4
e^{\pi i(\alpha+\beta)} \sin(\pi \alpha) \, \sin(\pi \beta) \,
\int_{-1}^1 f(t)\, dt\,.
\]

Changing the variable, $t=2x-1$, we get immediately the integral
defining the beta function, and  as a consequence
\begin{equation*} \label{I_1Main}
I_n(\alpha , \beta )= -2^{n+\alpha +\beta +3} \,  e^{\pi
i(\alpha+\beta)}\, \sin(\pi \alpha) \, \sin(\pi \beta) \,
\frac{\Gamma \left( \alpha +n+1\right) \Gamma \left( \beta
+n+1\right) }{\Gamma \left( 2n+ \alpha +\beta +2\right) } \,.
\end{equation*}
Using $\sin (\pi x) \Gamma(x) \Gamma(1-x)=\pi$, we obtain
(\ref{orthMain}) for $k=n$ and for $\alpha$ and $\beta$ with $\Re
\alpha $ and $\Re \beta $ positive. By analytic continuation the
identity (\ref{orthMain}) holds for every $\alpha, \beta \in
\mathbb C$.
\end{proof}

Observe that the right hand side in (\ref{orthMain}) vanishes for
$k=n$ if and only if either $-2n-\alpha-\beta-2$, or $n+\alpha$ or
$n+\beta$ is a non-negative integer. In some of these cases the
zero comes from integrating a single-valued and analytic function
along a curve in the region of analyticity; other values of
$\alpha$ and $\beta $ correspond to the special cases mentioned
before when there is a zero at $\pm 1$.

\section{A Riemann-Hilbert problem for Jacobi polynomials}
In this section we construct a Riemann-Hilbert problem whose
solution is given in terms of the Jacobi polynomials $P_n^{(\alpha,\beta)}$
with parameters satisfying
\begin{equation}\label{condition}
-n-\alpha-\beta \notin \N, \quad \text{and} \quad n+\alpha
    \notin \N, \quad \text{and} \quad n+\beta \notin \N.
\end{equation}

We consider $\Gamma$ as a curve in $\C$ with three points of
self-intersection. We let $\Gamma^o$ be the curve without the
points of self-intersection. The orientation of $\Gamma$, see Fig.
\ref{fig:pathGamma}, induces a $+$ and $-$ side in a neighborhood
of $\Gamma$, where the $+$ side is on the left while traversing
$\Gamma$ according to its orientation and the $-$ side is on the
right. We say that a function $Y$ on $\C \setminus \Gamma$ has a
boundary value $Y_+(t)$ for $t \in \Gamma^o$ if the limit of
$Y(z)$ as $z \to t$ with $z$ on the $+$ side of $\Gamma$ exist.
Similarly for $Y_-(t)$.

The Riemann-Hilbert problem asks for a $2\times 2$ matrix valued
function $Y : \C \setminus \Gamma \to \C^{2\times 2}$ such that
the following four conditions are satisfied.
\begin{enumerate}
\item[(a)] $Y$ is analytic on $\C \setminus \Gamma$. \item[(b)]
$Y$ has continuous boundary values on $\Gamma^o$, denoted by $Y_+$
and $Y_-$, such that
\[ Y_+(t) = Y_-(t)
    \begin{pmatrix} 1 & w(t; \alpha, \beta) \\
        0 & 1 \end{pmatrix}
        \qquad \mbox{for } t\in \Gamma^o. \]
\item[(c)] As $z\to\infty$,
\[ Y(z) = \left(I + O\left(\frac{1}{z}\right)\right)
    \begin{pmatrix} z^{n} & 0 \\ 0 & z^{-n} \end{pmatrix}. \]
\item[(d)]
$Y(z)$ remains bounded as $z \to t \in \Gamma \setminus \Gamma^o$.
\end{enumerate}
\medskip

This Riemann-Hilbert problem is similar to the Riemann-Hilbert
problem for orthogonal polynomials due to Fokas, Its, and Kitaev \cite{FIK},
see also \cite{Deiftbook}. Also the solution is similar.
It is built out of  the Jacobi polynomials $P_n^{(\alpha, \beta)}$ and
$P_{n-1}^{(\alpha, \beta)}$.
For parameters satisfying (\ref{condition}), the polynomial
$P_n^{(\alpha,\beta)}$ has degree $n$; recall that there is a degree reduction
if and only if
$-n-\alpha-\beta \in \{1,\ldots, n\}$. Therefore there is
a constant $c_n$ such that
\begin{equation} \label{defcn}
    c_n P_n^{(\alpha,\beta)} \quad \mbox{ is a monic polynomial}.
\end{equation}
From Theorem \ref{thm:main} and the condition (\ref{condition}) on
the parameters, it follows that $P_{n-1}^{(\alpha,\beta)}$
satisfies
\[ I_{n-1}(\alpha , \beta )=\int_{\Gamma} t^{n-1} P_{n-1}^{(\alpha,\beta)}(t) w(t; \alpha, \beta) dt
    \neq 0. \]
Thus we can define the constant $d_{n-1}=-2\pi i /I_{n-1}(\alpha ,
\beta )$ such that
\begin{equation} \label{defdn}
    d_{n-1} \int_{\Gamma} t^{n-1}  P_{n-1}^{(\alpha,\beta)}(t) w(t; \alpha, \beta) dt
    = - 2\pi i.
\end{equation}
Then we can state the following result.

\begin{proposition} \label{solRH}
The unique solution of the Riemann-Hilbert problem is given by
\begin{equation} \label{formulaY}
    Y(z) = \begin{pmatrix} c_n P_{n}^{(\alpha,\beta)}(z) &
    \frac{c_n}{2\pi i} \int_{\Gamma} \frac{P_{n}^{(\alpha,\beta)}(t)
        w(t; \alpha, \beta)}{t-z} dt \\[10pt]
    d_{n-1} P_{n-1}^{(\alpha,\beta)}(z) &
    \frac{d_{n-1}}{2\pi i} \int_{\Gamma}
    \frac{P_{n-1}^{(\alpha,\beta)}(t) w(t; \alpha, \beta)}{t-z} dt
    \end{pmatrix}.
\end{equation}
\end{proposition}
\begin{proof} The proof that (\ref{formulaY}) satisfies the
Riemann-Hilbert problem
is similar to the proof for usual orthogonal polynomials,
see e.g.\ \cite{Deiftbook,Kuijlaars}.
The condition (a) is obviously satisfied by (\ref{formulaY}).
The jump condition (b) follows from the Sokhotskii-Plemelj formula
\[ f_+(t) = f_-(t) + v(t), \qquad t \in \Gamma^o, \]
which is satisfied for $f(z) = \frac{1}{2\pi i} \int_{\Gamma}
\frac{v(t)}{t-z} dt$. The asymptotic condition (c) follows because
of the normalizations (\ref{defcn})--(\ref{defdn}) and the
orthogonality conditions given in Theorem \ref{thm:main}. For
example, for the $(2,2)$ entry of $Y$ the  condition (c) is
$Y_{22}(z) = z^{-n} + O(z^{-n-1})$ as $z \to \infty$, and this is
satisfied by the $(2,2)$ entry of (\ref{formulaY}) since
\begin{eqnarray*}
  \frac{d_{n-1}}{2\pi i} \int_{\Gamma}
    \frac{P_{n-1}^{(\alpha,\beta)}(t) w(t; \alpha, \beta)}{t-z} dt
    & = &
     - \sum_{k=0}^{\infty} \left(\frac{d_{n-1}}{2\pi i} \int_{\Gamma}
        t^k P_{n-1}^{(\alpha,\beta)}(t) w(t; \alpha,\beta) dt \right) z^{-k-1} \\
    & & =
    \left(
    - \frac{d_{n-1}}{2\pi i} \int_{\Gamma}
     t^{n-1} P_{n-1}^{(\alpha,\beta)}(t) w(t; \alpha,\beta) dt \right) z^{-n}
        + O(z^{-n-1}) \qquad \mbox{as } z \to \infty \\
    & & = z^{-n} + O(z^{-n-1}) \qquad \mbox{as } z \to \infty.
\end{eqnarray*}
For the first equality we used the Laurent expansion of
$\frac{1}{t-z}$ around $z=\infty$,  for the second
equality we used the orthogonality conditions satisfied
by $P_{n-1}^{(\alpha, \beta)}$ on $\Gamma$, and for the third equality
we used (\ref{defdn}).
Finally, the boundedness condition (d) is certainly satisfied for the
first column of (\ref{formulaY}). It is also satisfied by the second
column, since by analyticity we may deform the contour $\Gamma$ from
which it follows that the entries in the second column have analytic
continuations across $\Gamma$. Then they are certainly bounded.

The proof of uniqueness is as in \cite{Deiftbook,Kuijlaars}
and we omit it here.
\end{proof}

Using the Riemann-Hilbert problem, we can easily prove that the
orthogonality conditions in Theorem \ref{thm:main} characterize
the Jacobi polynomial in case the parameters satisfy
(\ref{condition}).
\begin{theorem} \label{thm:RH}
Assume that $\alpha, \beta  \in \C$ satisfy {\rm (\ref{condition})},
and that $\Gamma$, $w(z; \alpha ,\beta)$ are as above.
Then the monic Jacobi polynomial
$c_n P_n^{(\alpha,\beta)}$ is the only monic polynomial $p_n$ of degree $n$
that satisfies
\begin{equation} \label{orthopn}
    \int_{\Gamma} t^k p_n(t) w(t;\alpha, \beta) dt = 0,
    \qquad \mbox{ for } k=0,1, \ldots, n-1.
\end{equation}
\end{theorem}
\begin{proof}
The orthogonality conditions (\ref{orthopn}) are what is necessary
to fill the first row of $Y$. That is, if $p_n$ is a monic polynomial
of degree $n$, satisfying (\ref{orthopn}), then
\[  \begin{pmatrix} p_n(z) &
    \frac{1}{2\pi i} \int_{\Gamma} \frac{p_{n}(t)
        w(t; \alpha, \beta)}{t-z} dt \\[10pt]
    d_{n-1} P_{n-1}^{(\alpha,\beta)}(z) &
    \frac{d_{n-1}}{2\pi i} \int_{\Gamma}
    \frac{P_{n-1}^{(\alpha,\beta)}(t) w(t; \alpha, \beta)}{t-z} dt
    \end{pmatrix}
\]
satisfies all conditions in the Riemann-Hilbert problem for $Y$.
Since the solution is unique and given by (\ref{formulaY}) it follows
that $p_n(z) = Y_{11}(z) = c_n P_n^{(\alpha,\beta)}(z)$.
\end{proof}

\section{Non-hermitian quasiorthogonality}
In the rest of this paper we assume for simplicity that $\alpha$
and $\beta$ are real, although extension to non-real parameters is
possible. We also take $\alpha $, $\beta $, $\alpha +\beta  \notin
\Z$, so that (\ref{condition}) is automatically guaranteed.

Since in (\ref{orthMain}) we are integrating an analytic function,
we may deform the universal path $\Gamma$ freely within the region
of analyticity. In particular, if the integrand is integrable in the
neighborhood of a branch point ($\pm 1$ or $\infty$), we may allow
$\Gamma$ to pass through this point, taking care of using the
correct branch of the integrand.

%
%%%%%%%%%%%%%%%%%%%%%%%%%%%%%%%%%%%%%%%%%%%%%%%%%%%%%%%%%%
\begin{figure}[htb]
\centering \includegraphics[scale=0.5]{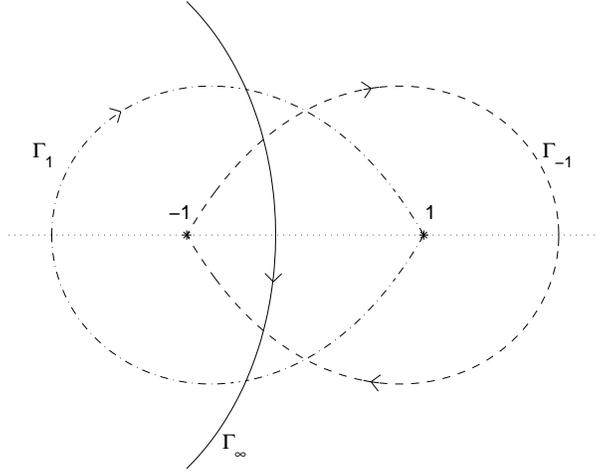} \caption{ Paths of
integration.}\label{fig:paths}
\end{figure}
%%%%%%%%%%%%%%%%%%%%%%%%%%%%%%%%%%%%%%%%%%%%%%%%%%%%%%%%%%
%

We define the following paths of integration (see
Fig.\ \ref{fig:paths}): $\Gamma_{1}$ will be an arbitrary curve
oriented clockwise, connecting $1 - i 0$ with $1 + i 0$ and lying
entirely in $\C \setminus [-1, +\infty)$, except for its
endpoints. The circle $\{ z \in \C:\, |z + 1|=2\}$, oriented
clockwise, is a good instance of a curve $\Gamma_{1}$.
Analogously, $\Gamma_{-1}$ will be an arbitrary curve oriented
clockwise, connecting $-1 + i 0$ with $-1 - i 0$ and lying
entirely in $\C \setminus (-\infty,1]$, except for its endpoints.
Finally, $\Gamma_{\infty}$  is a curve in $\C \setminus \big
((-\infty,-1] \cup [1, +\infty)\big)$, extending from $+ i \infty$
to $ - i \infty$ (for example, the imaginary axis, oriented
downward, might do as $\Gamma_{\infty}$).

In what follows, we denote by $\P_n$  the set of
all algebraic polynomials with complex coefficients of degree
$\leq n $ (for $n<0$ we assume $\P_n=\emptyset$), and by
$\P_n^*=\P_n \setminus \P_{n-1}$, the subset of polynomials of
degree exactly $n$. Also, for $c \in \mathbb R$, we use $[c]$
to denote the largest integer $\leq c$, as usual.

The following theorem shows that for certain combinations of
parameters a number of orthogonality relations are valid on
$\Gamma_{-1}$, $\Gamma_1$, or $\Gamma_{\infty}$. These
are called quasi-orthogonality relations, since in general
there are less than $n$ conditions on the polynomial
$P_n^{(\alpha,\beta)}$, so that these relations do not characterize
the polynomial.
When integrating over $\Gamma_{s}$, $s \in \{-1,1,\infty\}$, we mean by
$w(t; \alpha, \beta)$ a branch of the weight function $(1-t)^{\alpha}(t+1)^{\beta}$
which is continuous on $\Gamma_{s} \setminus \{s\}$.

\begin{theorem}
\label{thm:quasiorth}
 Assume that $\alpha, \beta, \alpha +\beta \in \R \setminus \Z$.
\begin{enumerate}
\item[\rm i)] If $n+\beta >-1$, and $k=\max\{0, [-\beta ] \}$, then
\begin{equation}\label{orthBeta}
    \int_{\Gamma_{-1}} q(t)\, P_n^{(\alpha , \beta )}(t) w(t; \alpha,
    k+\beta)\,
    dt \begin{cases}
=0, & \forall \, q \in \P_{n-k-1}\,,\\
\neq 0, & \forall \, q \in \P_{n-k}^*\,.
    \end{cases}
\end{equation}

\item[\rm ii)] If $n+\alpha >-1$, and $k=\max\{0, [-\alpha ] \}$, then
\begin{equation}\label{orthAlpha}
    \int_{\Gamma_{1}} q(t)\, P_n^{(\alpha , \beta )}(t) w(t; k+\alpha,
    \beta)\,
    dt \begin{cases}
=0, & \forall \, q \in \P_{n-k-1}\,,\\
\neq 0, & \forall \, q \in \P_{n-k}^*\,.
    \end{cases}
\end{equation}

\item[\rm iii)] If $n+\alpha +\beta <-1$ and $m=\min \{n-1,
[-(n+\alpha+ \beta+1)] \}$, then
\begin{equation}\label{orthInfty1}
 \int_{\Gamma_\infty} q(t)\, P_n^{(\alpha , \beta )}(t) w(t; \alpha,
    \beta)\, dt =0\,, \quad  \forall \, q \in \P_{m}\,.
\end{equation}
If additionally, $n+\alpha +\beta <-n-1$, then
\begin{equation*}\label{orthInfty2}
 \int_{\Gamma_\infty} q(t)\, P_n^{(\alpha , \beta )}(t) w(t; \alpha,
    \beta)\, dt
     \begin{cases}
     =0\,, & \forall \, q \in \P_{n-1}\,,  \\
     \neq 0\,, & \forall \, q \in \P_{n}^*\,.
     \end{cases}
\end{equation*}
\end{enumerate}
\end{theorem}
\begin{proof}
Assume that for a polynomial $q \in \P_n$ the function $ f(t)=q(t)
P_n^{(\alpha , \beta )}(t) w(t; \alpha, \beta)$ is integrable at
$t=-1$. Then we can deform the path $\Gamma$ in (\ref{orthMain})
into $\Gamma_{-1}$ traversed twice (in opposite directions), so
that
$$
\int_\Gamma f(t)\, dt=\left(e^{2\pi i \beta} -1\right)\,
\int_{\Gamma_{-1}} f(t)\, dt\,.
$$
If $n+\beta >-1$, then $k=\max\{0, [-\beta ] \} \in \{0, 1, \ldots,
n\}$, $k+\beta>-1$, and the integrability of $f$ at $-1$ is
guaranteed taking
$$
q(t)=(t+1)^k r(t)\,, \quad r \in \P_{n-k}\,.
$$
Then,
\begin{align*}
\int_\Gamma q(t) P_n^{(\alpha , \beta )}(t) w(t; \alpha, \beta)\,
dt &= \left(e^{2\pi i \beta} -1\right)\, \int_{\Gamma_{-1}} q(t)
P_n^{(\alpha , \beta )}(t) w(t; \alpha, \beta)\, dt\\
&=\left(e^{2\pi i \beta} -1\right)\, \int_{\Gamma_{-1}} r(t)
P_n^{(\alpha , \beta )}(t) w(t; \alpha, k+\beta)\, dt\,,
\end{align*}
and (\ref{orthBeta}) follows from Theorem \ref{thm:main}.

The proof of (\ref{orthAlpha}) follows by reversing the roles
of $\alpha$ and $\beta$.

The proof of part iii) uses similar arguments. If for a
polynomial $q $, $n+\alpha +\beta+ \deg q <-1$, then $ f=q
P_n^{(\alpha , \beta )} w(\cdot ; \alpha, \beta)$ is integrable at
$\infty$, and we can deform the path $\Gamma$ in (\ref{orthMain})
into $\Gamma_{\infty}$ traversed four times (two in each
direction), each time with a different branch of the integrand.
Then again
$$
\int_{\Gamma} q(t)\, P_n^{(\alpha , \beta )}(t) w(t; \alpha,
    \beta)\, dt = c\, \int_{\Gamma_\infty} q(t)\, P_n^{(\alpha , \beta )}(t) w(t;
\alpha,
    \beta)\, dt\,, \quad c \neq 0.
$$
In particular, $n+\alpha +\beta+ m <-1$, and we may apply
(\ref{orthMain}) with $k=0, \ldots, m \leq n-1$ to establish
(\ref{orthInfty1}). Furthermore, if $n+\alpha +\beta <-n-1$, then
(\ref{orthMain}) can be used up to $k=n$. This concludes the
proof.
\end{proof}

\medskip

\textsc{Remarks:} If $-1 < n+\beta <0$, then $\max\{0, [-\beta ]
\}=n$, and (\ref{orthBeta}) is reduced to a single condition,
$$
 \int_{\Gamma_{-1}}  P_n^{(\alpha , \beta )}(t) w(t; \alpha,
    n+\beta)\,
    dt \neq 0\,.
$$
An analogous degenerate situation is observed when
$-1 < n+\alpha <0$. On the other hand, if $-n-1< n+\alpha +\beta <-1$, the
integral in (\ref{orthInfty1}) diverges for $q \in \P_{m+1}^*$.

\section{Orthogonality on a single contour}

Sometimes it is possible to obtain a full set of orthogonality
conditions on $\Gamma_{-1}$, $\Gamma_1$ or $\Gamma_{\infty}$.

\begin{theorem}[Non-hermitian orthogonality] \label{thm:nonherm}
If for $\alpha , \beta, \alpha +\beta  \in \R\setminus \Z$, at
least one of the following conditions is fulfilled:
\begin{equation}\label{3conditions}
\alpha >-1, \qquad  \beta >-1,  \qquad 2n+ \alpha +\beta< 0,
\end{equation}
then Jacobi polynomials $P_n^{(\alpha , \beta )}$ satisfy a full
set of non-hermitian (complex) orthogonality conditions:
\begin{equation}\label{orthNonHerm}
    \int_{\gamma} q(t)\, P_n^{(\alpha , \beta )}(t) w(t; \alpha,
    \beta)\,
    dt \begin{cases}
    =0\,, & q \in \P_{n-1},\\
    \neq 0 \,, & q \in \P_{n}^*,
    \end{cases}
\end{equation}
where
$$
  \gamma= \begin{cases}
\Gamma_{1}\,, & \text{if }\; \alpha  >-1,\\
\Gamma_{-1} \,, & \text{if } \; \beta >-1, \\
\Gamma_{\infty} \,, & \text{if } \; 2n+ \alpha +\beta <0 \,.
    \end{cases}
$$
(If $-1<2n+\alpha +\beta <0$ the integral in {\rm (\ref{orthNonHerm})}
diverges for $\gamma=\Gamma_{\infty}$ and $q \in \P^*_{n}$).

The conditions {\rm (\ref{orthNonHerm})} characterize the Jacobi polynomial
$P_n^{(\alpha , \beta )}$ of degree $n$ up to a constant factor.
\end{theorem}
\begin{proof}
This is an immediate consequence of Theorem \ref{thm:quasiorth},
since under our assumptions, $k=0$ in
(\ref{orthBeta})--(\ref{orthAlpha}). It is straightforward to show
that orthogonality conditions in (\ref{orthNonHerm}) are
equivalent to (\ref{orthopn}), and thus by Theorem \ref{thm:RH}
characterize the polynomial.
\end{proof}

\medskip

The following result, describing the real orthogonality of Jacobi
polynomials, is classical, but we put it within the general
framework.
\begin{corollary}[Real orthogonality] \label{thm:real}
If for $\alpha , \beta, \alpha +\beta  \in \R\setminus \Z$, two of
the three conditions in {\rm (\ref{3conditions})} are fulfilled, then
Jacobi polynomials $P_n^{(\alpha , \beta )}$ satisfy a full set of
orthogonality conditions on the real line:
\begin{equation}\label{orth1C1}
    \int_{\Delta} q(t)\, P_n^{(\alpha , \beta )}(t) w(t; \alpha,
    \beta)\,
    dt \begin{cases}
=0\,, & q \in \P_{n-1},\\
\neq 0 \,, & q \in \P_{n}^*,
    \end{cases}
\end{equation}
where
\begin{equation} \label{delta}
  \begin{cases}
[-1,1]\,, & \text{if }\; \alpha , \beta >-1,\\
(-\infty, -1] \,, & \text{if } \; 2n+ \alpha +\beta <0 \text{ and } \beta >-1, \\
\, [1, +\infty) \,, & \text{if } \; 2n+ \alpha +\beta <0 \text{
and } \alpha  >-1\,.
    \end{cases}
\end{equation}
(If $-1<2n+\alpha +\beta <0$ the integral in {\rm (\ref{orth1C1})}
diverges for $q \in \P^*_{n}$).

The conditions {\rm (\ref{orth1C1})} characterize the Jacobi polynomial
$P_n^{(\alpha , \beta )}$ of degree $n$ up to a constant factor.
\end{corollary}
\begin{proof}
This is a consequence of Theorem \ref{thm:nonherm}. For instance, if
$2n+\alpha + \beta < 0$, by (\ref{orthNonHerm}) we have non-hermitian
orthogonality on $\Gamma_{\infty}$; if additionally $\alpha>-1$, we
can deform $\Gamma_{\infty}$ into $[1,+\infty)$ traversed twice, and the
statement follows.
\end{proof}

\section{Multiple orthogonality}

Theorem \ref{thm:nonherm} provides orthogonality conditions,
characterizing Jacobi polynomials when their parameters belong to
the region in the $(\alpha , \beta )$-plane given by at least
one of the conditions (\ref{3conditions}).

For other combinations of parameters we still have some
orthogonality relations according to Theorem \ref{thm:quasiorth},
but each of the three orthogonality relations (\ref{orthBeta}),
(\ref{orthAlpha}), (\ref{orthInfty1}) does not give enough
conditions to determine $P_n^{(\alpha,\beta)}$ by itself.
However, the three relations taken together give
$n$ or more relations for $P_n^{(\alpha,\beta)}$.
They constitute what we call a set of multiple orthogonality
conditions.
In many cases there will be more than $n$ conditions,
so that the relations of Theorem \ref{thm:quasiorth}
overdetermine $P_n^{(\alpha,\beta)}$.

We are going to discuss this in more detail now.

\subsection{Multiple orthogonality as an alternative to orthogonality on
a single contour}

We will consider the following subcases of the situation,
described in Theorem \ref{thm:nonherm}.

\begin{theorem} \label{prop51}
Let $\alpha, \beta, \alpha +\beta  \in \R\setminus \Z$ such that
exactly one of the conditions {\rm (\ref{3conditions})}
is satisfied, and such that either
$-n < \alpha < -1$ or $-n < \beta < -1$. Then
\begin{gather}
    \int_{\gamma} q(t)\,
P_n^{(\alpha , \beta)}(t) w(t; \alpha ,
    \beta)\,
    dt =0\,, \quad q \in \P_{k-1} \,,  \label{11}\\
\intertext{and }
 \int_{\Delta} q(t)\, P_n^{(\alpha , \beta )}(t)
w(t; \alpha_1, \beta_1)\, dt \, \begin{cases}
=0\,, & q \in \P_{n-k-1},\\
\neq 0 \,, & q \in \P_{n-k}^*,
    \end{cases}
\label{12}
\end{gather}
where the corresponding parameters are gathered in the following
table:
%\begin{center}
\begin{equation*}\label{table2}
\begin{tabular}{|c|c|c|c|c|c|}
  \hline
   \strut Cases& $\gamma $ & $\Delta$ & $\alpha_1$ & $\beta_1$ &  k \\
  \hline
  $\alpha >-1$, $-n<\beta <-1$ \strut &  $\Gamma_{1}$ & $[-1,1]$ & $\alpha $ & $\beta +[-\beta]$ & $[-\beta]$ \\
  $\beta >-1$, $-n<\alpha <-1$ \strut &  $\Gamma_{-1}$ & $[-1,1]$ & $\alpha+[-\alpha] $ & $\beta $ & $[-\alpha]$ \\
 $2n+\alpha+\beta <0$, $-n<\beta <-1$\strut &  $\Gamma_{\infty}$ & $(-\infty,-1]$ & $\alpha $ & $\beta +[-\beta]$ & $[-\beta]$ \\
  $2n+\alpha + \beta <0$, $-n<\alpha <-1$ \strut &  $\Gamma_{\infty}$ & $[1,+\infty)$ & $\alpha+[-\alpha] $ & $\beta $ & $[-\alpha]$ \\
 \hline
\end{tabular}
\end{equation*}
(If $-1<2n+\alpha +\beta <0$ the integral in {\rm (\ref{12})} diverges
for $q \in \P^*_{n-k}$).

In each case, these multiple orthogonality conditions characterize
the Jacobi polynomial $P_n^{(\alpha , \beta )}$ of degree $n$ up
to a constant factor.
\end{theorem}
\begin{proof}
Assume for instance that $\beta >-1$ and $-n < \alpha <-1$.
Then we have on the one hand by Theorem \ref{thm:nonherm} that
\begin{equation}\label{Prop51_f1}
    \int_{\Gamma_{-1}} q(t)\, P_n^{(\alpha , \beta )}(t) w(t; \alpha,
    \beta)\,
    dt \begin{cases}
=0\,, & q \in \P_{n-1},\\
\neq 0 \,, & q \in \P_{n}^*\,,
    \end{cases}
\end{equation}
and on the other hand by Theorem \ref{thm:quasiorth} ii), that
\begin{equation}\label{Prop51_f2}
    \int_{\Gamma_{1}} q(t)\, P_n^{(\alpha , \beta )}(t) w(t;  \alpha+k,
    \beta)\,
    dt \begin{cases}
=0\,, & q \in \P_{n-k-1},\\
\neq 0 \,, & q \in \P_{n-k}^*,
    \end{cases}   %\quad \text{with } \Sigma \in \F_{1}\,,
\end{equation}
where $k = [-\alpha]$.
By Euclid's algorithm, any $q \in \P_m$, $m \geq k$, may be
represented in the form $q(t)=(t-1)^{k} r(t)+s(t)$ with $r \in
\P_{m-k}$ and $s \in \P_{k-1}$. Thus, (\ref{Prop51_f1}) is equivalent
to (\ref{11}) and
\begin{equation}\label{Prop51_f3}
    \int_{\Gamma_{-1}} q(t)\, P_n^{(\alpha , \beta )}(t) w(t; \alpha+k,
    \beta)\,
    dt \begin{cases}
=0\,, & q \in \P_{n-k-1},\\
\neq 0 \,, & q \in \P_{n-k}^*.
    \end{cases}
\end{equation}
Since in (\ref{Prop51_f2})--(\ref{Prop51_f3}) we can deform the
path of integration and the functions are integrable at the
singularities $\pm 1$, these two identities are equivalent to (\ref{12})
with $\Delta = [-1,1]$, $\alpha_1 = \alpha + [-\alpha]$, and
$\beta_1 = \beta$.
Thus, (\ref{11})--(\ref{12}) are equivalent to (\ref{Prop51_f1}),
and they characterize $P_n^{(\alpha , \beta )}$ up to a constant
factor by Theorem \ref{thm:nonherm}.

The other cases are handled in a similar fashion.
\end{proof}

\medskip

\begin{theorem} \label{prop52}
Let  $\alpha , \beta, \alpha +\beta  \in \R\setminus \Z$ be such
that $-n<\alpha+\beta +n <-1$ and either $\alpha >-1$ or $\beta
>-1$. Then with $m=[-(n+\alpha+ \beta+1)]$,
\begin{gather*}
    \int_{\Delta} q(t)\,
P_n^{(\alpha , \beta)}(t) w(t; \alpha ,
    \beta)\,
    dt =0\,, \quad q \in \P_{m} \,,  \label{21}\\
\intertext{and }
 \int_{\gamma} t^k\, P_n^{(\alpha , \beta )}(t)
w(t; \alpha, \beta)\, dt \, \begin{cases}
=0\,, & \text{if } k \in \N, \; m+1 \leq k \leq n-1\,,\\
\neq 0 \,, & \text{if } k=n\,,
    \end{cases}
\label{22} \\ \intertext{where} \begin{cases} \Delta=[1,+\infty),
\; \gamma=\Gamma_{1}, & \text{if } \alpha >-1, \\
\Delta=(-\infty,-1], \; \gamma=\Gamma_{-1}, & \text{if } \beta
>-1.
\end{cases} \nonumber
\end{gather*}
In each case, these multiple orthogonality conditions characterize
the Jacobi polynomial $P_n^{(\alpha , \beta )}$ of degree $n$ up
to a constant factor.
\end{theorem}
\begin{proof}
This is a corollary of Theorem \ref{thm:nonherm} and
(\ref{orthInfty1}).
\end{proof}

\subsection{Multiple orthogonality when there is no orthogonality on a single contour}

As we have seen, in the cases analyzed so far the multiple
orthogonality was in a certain sense ``optional'': whenever at least one of
the conditions in (\ref{3conditions}) is satisfied, we can restrict
ourselves to orthogonality on a single contour $\Gamma_{-1}$, $\Gamma_1$
or $\Gamma_{\infty}$, as in Theorem \ref{thm:nonherm}. In the remaining cases
there are quasi-orthogonality conditions on $\Gamma_{-1}$, $\Gamma_1$
and $\Gamma_{\infty}$ as specified in Theorem \ref{thm:quasiorth}, but we
need at least two of these sets of quasi-orthogonality relations to
characterize the Jacobi polynomial.

So in this part we assume that $\alpha < -1$, $\beta < -1$ and $2n + \alpha + \beta > 0$.
In our next result we consider cases where the parameters are such that
a combination of two of the cases in Theorem \ref{thm:quasiorth} give at least
$n$ orthogonality conditions. The theorem says how to obtain from that
$n$ conditions that characterize $P_n^{(\alpha,\beta)}$.

\begin{theorem} \label{prop54}
Let $\alpha , \beta, \alpha +\beta  \in \R\setminus \Z$ such that
$\alpha < -1$, $\beta < -1$ and $2n + \alpha + \beta > 0$.
\begin{enumerate}
\item[\rm i)] If $\alpha+\beta +n >-1$, then
\begin{gather}
    \int_{\Gamma_{-1}} q(t)\,
P_n^{(\alpha , \beta)}(t) w(t; \alpha ,
    \beta+[-\beta ])\,
    dt =0\,, \quad q \in \P_{[-\alpha ]-1} \,,  \label{41}\\
    \int_{\Gamma_{1}} q(t)\,
P_n^{(\alpha , \beta)}(t) w(t; \alpha +[-\alpha ],
    \beta)\,
    dt =0\,, \quad q \in \P_{[-\beta]-1} \,,  \label{42} \\
\intertext{and }
 \int_{-1}^1 q(t)\, P_n^{(\alpha , \beta )}(t)
w(t; \alpha+[-\alpha ], \beta+ [-\beta])\, dt \, \begin{cases}
=0\,, & q \in \P_{n-[-\alpha]-[-\beta]-1}\,,\\
\neq 0 \,, & q \in \P_{n-[-\alpha]-[-\beta]}\,.
    \end{cases}
\label{43}
\end{gather}
\item[\rm ii)] If $\alpha < -n$, then
    \begin{gather}
    \int_{-\infty}^{-1} q(t)\,
P_n^{(\alpha , \beta)}(t) w(t; \alpha, \beta+[-\beta ])\,
    dt =0\,, \quad q \in \P_{[-\alpha]-n-1} \,,  \label{tm41}\\
    \int_{\Gamma_{\infty}} q(t)\,
P_n^{(\alpha , \beta)}(t) w(t; \alpha ,
    \beta)\,
    dt =0\,, \quad q \in \P_{[-\beta]-1} \,,  \label{tm42} \\
\intertext{and }
 \int_{\Gamma_{-1}} q(t) \, P_n^{(\alpha , \beta )}(t)
 t^{-n+[-\alpha]} w(t; \alpha, \beta+ [-\beta])\, dt \, \begin{cases}
=0\,, &  q \in \P_{2n-[-\alpha]-[-\beta]-1}\,, \\
\neq 0 \,, & q \in \P_{2n-[-\alpha]-[-\beta]}\,.
    \end{cases}
\label{tm43}
\end{gather}
    \item[\rm iii)] If $\beta < -n$, then
    \begin{gather}
    \int_{1}^{+\infty} q(t)\,
P_n^{(\alpha , \beta)}(t) w(t; \alpha+[-\alpha], \beta)\,
    dt =0\,, \quad q \in \P_{[-\beta]-n-1} \,,  \label{tm44}\\
    \int_{\Gamma_{\infty}} q(t)\,
P_n^{(\alpha , \beta)}(t) w(t; \alpha, \beta)\,
    dt =0\,, \quad q \in \P_{[-\alpha]-1} \,,  \label{tm45} \\
\intertext{and }
 \int_{\Gamma_{1}} q(t) \, P_n^{(\alpha , \beta )}(t) t^{-n+[-\beta]}
w(t; \alpha+[-\alpha], \beta)\, dt \, \begin{cases}
=0\,, &  q \in \P_{2n-[-\alpha]-[-\beta]-1}\,,  \\
\neq 0 \,, & q \in \P_{2n-[-\alpha]-[-\beta]}\,.
    \end{cases}
\label{tm46}
\end{gather}
\end{enumerate}

In all three cases, these orthogonality conditions characterize the Jacobi polynomial
$P_n^{(\alpha , \beta )}$ of degree $n$ up to a constant factor.
\end{theorem}
\begin{proof}
Indeed, by parts i) and ii) of Theorem \ref{thm:quasiorth} we have
\begin{gather}
    \int_{\Gamma_{-1}} q(t)\, P_n^{(\alpha , \beta )}(t) w(t; \alpha,
    \beta+[-\beta ])\,
    dt \begin{cases}
=0, & \, q \in \P_{n-[-\beta ]-1}\,,\\
\neq 0, &  \, q \in \P_{n-[-\beta ]}^*\,,
    \end{cases} \label{thm63_1}
\\ \intertext{and}
    \int_{\Gamma_{1}} q(t)\, P_n^{(\alpha , \beta )}(t) w(t; \alpha+[-\alpha ],
    \beta)\,
    dt \begin{cases}
=0, &  \, q \in \P_{n-[-\alpha ]-1}\,,\\
\neq 0, &  \, q \in \P_{n-[-\alpha ]}^*\,.
    \end{cases} \label{thm63_2}
\end{gather}
Since $[-\alpha] + [-\beta] \leq n$, we obtain (\ref{41}) and (\ref{42}).
For $q \in \P_{n-[-\alpha]-[-\beta]}$, we take the polynomial $(t-1)^{[-\alpha]} q(t)$
in (\ref{41}) and deform the contour $\Gamma_{-1}$ to $[-1,1]$, to obtain (\ref{43}) as well.

To prove that the conditions characterize the Jacobi polynomial, let $p_n$ be a polynomial
of degree $\leq n$ that satisfies (\ref{41}), (\ref{42}), (\ref{43}).
These conditions imply (\ref{thm63_1}) and (\ref{thm63_2}), and passing to the contour
$\Gamma$, we obtain
\[ \int_{\Gamma} q(t) (t+1)^{[-\beta]} p_n(t) w(t,\alpha, \beta) dt = 0,
    \quad  q \in \P_{n-[-\beta]-1}, \]
and
\[ \int_{\Gamma} q(t) (t-1)^{[-\alpha]} p_n(t) w(t,\alpha,\beta) dt = 0,
    \quad q \in \P_{n-[-\alpha]-1}. \]
Thus
\begin{equation} \label{thm63_3}
\int_{\Gamma} q(t) p_n(t) w(t,\alpha,\beta) dt = 0
\end{equation}
for every $q \in \P_{n-1}$ of the form $q(t)=(t+1)^{[-\beta]} r(t)+ (t-1)^{[-\alpha]} s(t)$
with $r \in \P_{n-[-\beta]-1}$ and $s \in \P_{n-[-\alpha]-1}$.
The linear mapping
\begin{equation} \label{mapping}
    \P_{n-[-\beta]-1} \times \P_{n-[-\alpha]-1} \to \P_{n-1} :
    (r,s) \mapsto q(t) = (t+1)^{[-\beta]} r(t) + (t-1)^{[-\alpha]} s(t)
\end{equation}
has a kernel consisting of pairs $(r,s)$ such that $(t+1)^{[-\beta]} r(t) + (t-1)^{[-\alpha]} s(t) \equiv 0$.
Then it is easy to see that $r(t) = (t-1)^{[-\alpha]} p(t)$
and $s(t) = -(t+1)^{[-\beta]} p(t)$ with $p \in \P_{n-[-\alpha]-[-\beta]-1}$.
So the kernel of the mapping (\ref{mapping}) has dimension $n-[-\alpha]-[-\beta]$.
Then by the dimension theorem from linear algebra the range of the mapping (\ref{mapping})
has dimension $(n-[-\beta]) + (n-[-\alpha]) - (n-[-\alpha]-[-\beta]) = n$,
and so the mapping is surjective.
It follows that (\ref{thm63_3}) holds for every $q \in \P_{n-1}$,
and then it follows from Theorem \ref{thm:RH} that $p_n$ is a multiple of the Jacobi
polynomial of degree $n$, so that the orthogonality conditions (\ref{41}), (\ref{42}), (\ref{43})
characterize $P_n^{(\alpha,\beta)}$ up to a constant factor.
This proves part i) of the theorem.

The other parts are proved in a similar way.
\end{proof}
\medskip

What remains is the case where the parameters are such that we need all
three cases of Theorem \ref{thm:quasiorth} in order to obtain at least $n$
conditions on $P_n^{(\alpha,\beta)}$.

\begin{theorem} \label{thm:last}
Let $\alpha , \beta, \alpha +\beta  \in \R\setminus \Z$ be such that
$\alpha > -n$, $\beta > -n$ and $\alpha+\beta +n <-1$.
Then the following hold:
\begin{gather}
    \int_{\Gamma_{-1}}q(t) \, P_n^{(\alpha , \beta )}(t)
        w(t; \alpha, \beta+[-\beta ])\, dt \, =0\,,
            \quad q \in \P_{n-[-\beta]-1} \,,  \label{tm47}\\
    \int_{\Gamma_{1}} q(t) \, P_n^{(\alpha , \beta )}(t)
        w(t; \alpha+[-\alpha], \beta)\, dt \, =0\,,
            \quad q \in \P_{n-[-\alpha]-1} \,, \label{tm48} \\
    \intertext{and }
    \int_{\Gamma_{\infty}} q(t) \, P_n^{(\alpha , \beta )}(t)
        w(t; \alpha, \beta)\, dt \, =0\,,
            \qquad q \in \P_{[-\alpha]+[-\beta]-n-1}\,, \label{inftyLast}
\end{gather}
and these  orthogonality conditions characterize the
Jacobi polynomial $P_n^{(\alpha , \beta )}$ of degree $n$ up to a
constant factor.
\end{theorem}

\begin{proof}
The relations (\ref{tm47})--(\ref{inftyLast}) were already established
in (\ref{orthBeta})--(\ref{orthInfty1}) from Theorem
\ref{thm:quasiorth}.

To show that these conditions characterize the Jacobi polynomial,
we assume that $p_n$ is a polynomial of degree $\leq n$, that
satisfies (\ref{tm47})--(\ref{inftyLast}). Passing to the
contour $\Gamma$, we then get that
\begin{equation} \label{qorthogonality}
    \int_{\Gamma} q(t) p_n(t) w(t;\alpha, \beta) dt = 0
\end{equation}
for every polynomial $q$ satisfying either
$q(t) = (t+1)^{[-\beta]} r(t)$ with $r \in \P_{n-[-\beta]-1}$, or
$q(t) = (t-1)^{[-\alpha]} s(t)$ with $ s \in \P_{n-[-\alpha]-1}$, or
$q(t) \in \P_{[-\alpha]+[-\beta]-n-1}$.

In the rest of the proof we will show that every polynomial
of degree $\leq n-1$ can be written in the form
\[ (t+1)^{[-\beta]} r(t) + (t-1)^{[-\alpha]} s(t) + q(t) \]
with $r \in \P_{n-[-\beta]-1}$, $s \in \P_{n-[-\alpha]-1}$ and $q \in \P_{[-\alpha]+[-\beta]-n-1]}$.
Having that, we get that (\ref{qorthogonality}) holds for
every $q \in \P_{n-1}$, and this characterizes the Jacobi polynomial
up to a constant factor by Theorem \ref{thm:RH}.

So we want to prove that the linear mapping
\begin{equation} \label{mapping2}
\begin{array}{l}
\P_{n-[-\beta]-1} \times \P_{n-[-\alpha]-1} \times \P_{[-\alpha]+[-\beta]-n-1} \to
\P_{n-1} : \\[10pt]
(r,s,q) \mapsto (t+1)^{[-\beta]} r(t) + (t-1)^{[-\alpha]} s(t) + q(t)
\end{array}
\end{equation}
is surjective. Since it is a mapping between $n$-dimensional
vector spaces, it suffices to show that (\ref{mapping2}) is injective.

Assume that $(r,s,q)$ belongs to the kernel of (\ref{mapping2}).
We are going to count the possible sign changes among the coefficients
of $(t+1)^{[-\beta]}r(t) +q(t) = - (t-1)^{[-\alpha]}s(t)$ with respect
to the standard monomial basis.
Here we use the following lemma which can be found, for instance, in
\cite[Section V, Ch.\ 1, problems 4 and 32]{polya/szego:1954}.
\begin{lemma} \label{polya} If
$$
p(z) = \sum_{j=0}^n a_j z^j \quad \text{and} \quad (z+1) p(z) =
\sum_{j=0}^{n+1} b_j z^j\,,
$$
then the number of sign changes among the $b_j$'s is not more than
the number of sign changes among the $a_j$'s.
\end{lemma}

Applying the lemma $[-\beta]$ times, we get that the number of sign changes
among the coefficients of $(t+1)^{[-\beta]} r(t)$ is at most the
number of sign changes among the coefficients of $r(t)$, which is at most
$n-[-\beta]-1$, since $r \in \P_{n-[-\beta]-1}$. When we add a polynomial
of degree $\leq k-1$, the number of sign changes among the coefficients
can increase with at most $k$.
Hence, $-(t-1)^{[-\alpha]} s(t)$ has at most
$n-[-\beta]-1 + [-\alpha]+[-\beta]-n = [-\alpha] - 1$ sign changes among
its coefficients.
Then by the Descartes rule of signs, see e.g.\ \cite[Section V, Ch.~1, problem 36]{polya/szego:1954},
$-(t-1)^{[-\alpha]} s(t)$ has at most $[-\alpha]-1$ positive real
zeros, counted according to their multiplicities, unless it is identically zero.
So we conclude that $s \equiv 0$. Then also $r \equiv 0$ and $q \equiv 0$.
This shows that (\ref{mapping2}) is indeed injective, which completes the
proof of Theorem \ref{thm:last}.
\end{proof}

%%%%%%%%%%%%%%%%%%%%%%%%%%%%%%%%%%%%%%%%%%%%%%%%%%%%%%%%%%
\begin{figure}[htb]
\centering \includegraphics[scale=0.9]{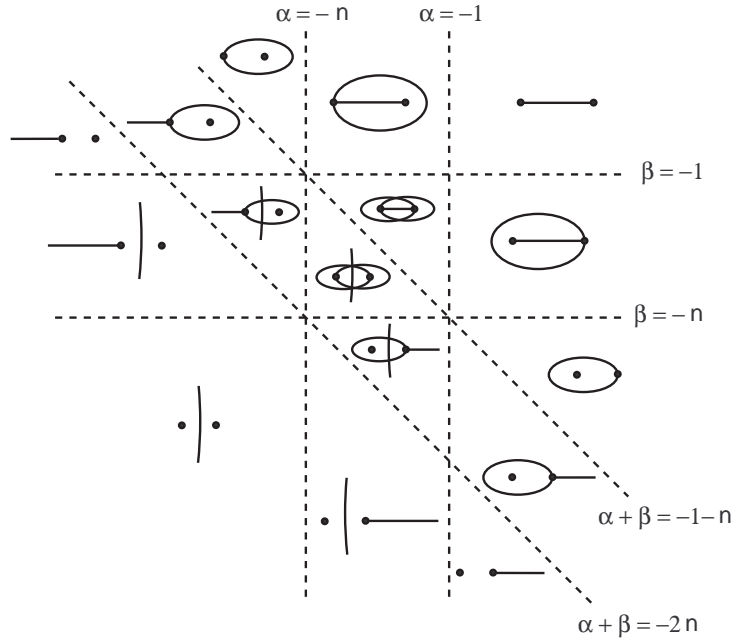}
% \unitlength=1in \psfragscanon
% \psfrag*{aka}{$\alpha =-n$}
% \psfrag*{bb}{$\alpha =-1$}
  \caption{Contours with orthogonality
conditions (black dots denote $-1$ and
$1$).}\label{fig:conditions}
\end{figure}
%%%%%%%%%%%%%%%%%%%%%%%%%%%%%%%%%%%%%%%%%%%%%%%%%%%%%%%%%%
%

To summarize, the paths of orthogonality, corresponding to
different cases studied, are represented schematically in
Fig.~\ref{fig:conditions}.

\section{Zeros}

Finally, it is interesting to discuss the implication of the
orthogonality relations derived to the location of the zeros of
the Jacobi polynomials.

Standard (hermitian) orthogonality conditions with respect to a
positive measure on the real line yield a lower bound (and in some
situation, the exact number) of different zeros of the orthogonal
polynomial on the convex hull of the support of the measure of
orthogonality.
For instance, by a classical argument, if two of the three
conditions in (\ref{3conditions}) are fulfilled, orthogonality
(\ref{orth1C1}) ensures that $P_n^{(\alpha , \beta )}$ has exactly
$n$ real and simple zeros, all located in the interval $\Delta$
given by (\ref{delta}).

In the situation of Theorem \ref{prop51}, that is when $(\alpha ,
\beta )$ or $(\beta , \alpha ) \in (-1, +\infty) \times (-n,-1)$,
quasi-orthogonality (\ref{12}) implies that $ P_n^{(\alpha , \beta
)}$ has \emph{at least} $n-k$ different zeros on $(-1,1)$, where
$k=[-\alpha ]$ if $-n<\alpha <-1$, and $k=[-\beta ]$ if $ -n<\beta
<-1$, see \cite{Brezinski02}.

Finally, in the situation of Theorem \ref{prop54}, $\alpha, \beta
\in (-n, -1)$ and $\alpha+\beta +n >-1$, we have that $
P_n^{(\alpha , \beta )}$ has at least $n-[-\alpha]-[-\beta]$
different zeros on $(-1,1)$.

In the other situations we have no quasi-orthogonality on
the real line, and the information on the real zeros is void.

It is interesting to compare these values with the well-known
Hilbert-Klein formulas \cite[Theorem 6.72]{szego:1975} that give
us the exact number of zeros of $P_n^{(\alpha ,\beta) }$ on the
real line. It is sufficient to restrict our attention to $[-1,1]$,
since the other two intervals of $\R$ can be analyzed similarly.

Following \cite{szego:1975}, let
$$
E(u)= \begin{cases}
 0, & \text{if } u \leq 0, \\ \,
 u-1, & \text{if } u \in \N, \\ \,
 [u], & \text{otherwise,}
\end{cases}
$$
and let $\even(u)$ (resp., $\odd(u)$) be the even (resp., odd)
value from the set $\{ u, u+1\}$, when $u \in \N$. Then the number
of zeros of $P_n^{(\alpha , \beta ))}$ in $(-1,1)$ is given by
\begin{equation}\label{N}
N=\begin{cases} \even\left(E\left( \frac{|2n+\alpha +\beta
+1|-|\alpha |
-|\beta|+1 }{2}\right)\right), & \text{if } \kappa_n(\alpha , \beta )>0, \\
\odd\left(E\left( \frac{|2n+\alpha +\beta +1|-|\alpha | -|\beta|+1
}{2}\right)\right), & \text{if } \kappa_n(\alpha , \beta )<0,
\end{cases}
\end{equation}
where $\kappa_n(\alpha , \beta )=(-1)^n (\alpha +1)\cdots(\alpha
+n)(\beta +1)\cdots (\beta +n)$.

It is straightforward to check that $N$ is positive only in one of the
following cases:
\begin{itemize}
\item $\alpha, \beta >-1$. Then $N=n$, as ensured by
(\ref{orth1C1}).
 \item $(\alpha , \beta )$ or $(\beta, \alpha ) \in (-1, +\infty) \times
(-n,-1)$. For instance, if $\beta >-1$ and $-n<\alpha <-1$, then
$$
E\left( \frac{|2n+\alpha +\beta +1|-|\alpha | -|\beta|+1
}{2}\right)=[n+\alpha +1]=n-[-\alpha ]\,,
$$
since for $\alpha <0$, $[\alpha ]+[-\alpha ]=-1$. In other words,
this value matches the lower bound on the number of zeros $N$,
predicted by the quasi-orthogonality relation on $[-1,1]$ given
in Theorem \ref{prop51}.

\item If $\alpha, \beta \in (-n, -1)$ and $\alpha+\beta +n >-1$,
$$ E\left( \frac{|2n+\alpha +\beta
+1|-|\alpha | -|\beta|+1 }{2}\right)=\begin{cases} [n+\alpha
+\beta +1]=n-[-\alpha-\beta  ]\,, & \text{if } \alpha +\beta
\notin \Z\,, \\
n+\alpha +\beta , & \text{otherwise.}
\end{cases}
$$
Observe that for $\tilde \alpha=\alpha +[-\alpha ]$, $\tilde \beta
=\beta  +[-\beta  ]$,
\begin{equation*}
\begin{split}
 &(-1)^{n-[-\alpha ]-[-\beta ]}\kappa_n(\alpha ,\beta ) \\
 &= (-\alpha
 -1)\cdots (-\tilde \alpha) (\tilde \alpha +1)\cdots (\alpha +n) (-\beta
 -1)\cdots (-\tilde \beta) (\tilde \beta +1)\cdots (\beta +n)>0\,.
\end{split}
\end{equation*}
Since $n-[-\alpha-\beta  ] \leq n-[-\alpha]-[-\beta  ] \leq
n-[-\alpha-\beta  ]+1$, we get from (\ref{N}) that $N=
n-[-\alpha]-[-\beta  ]$, which coincides again with the lower
bound on the number of zeros, predicted by Theorem \ref{prop54}.
\end{itemize}

Thus, quasi-orthogonality sheds a new light on the Hilbert-Klein
formulas, explaining the lower bound on the number of zeros on
$\R$.

%%%%%%%%%%%%%%%%%%%%%%%%%%%%%%%%%%%%%%%%%%%%%%%%%%%%%%%%%%
\begin{figure}[htb]
\centering \includegraphics[scale=0.7]{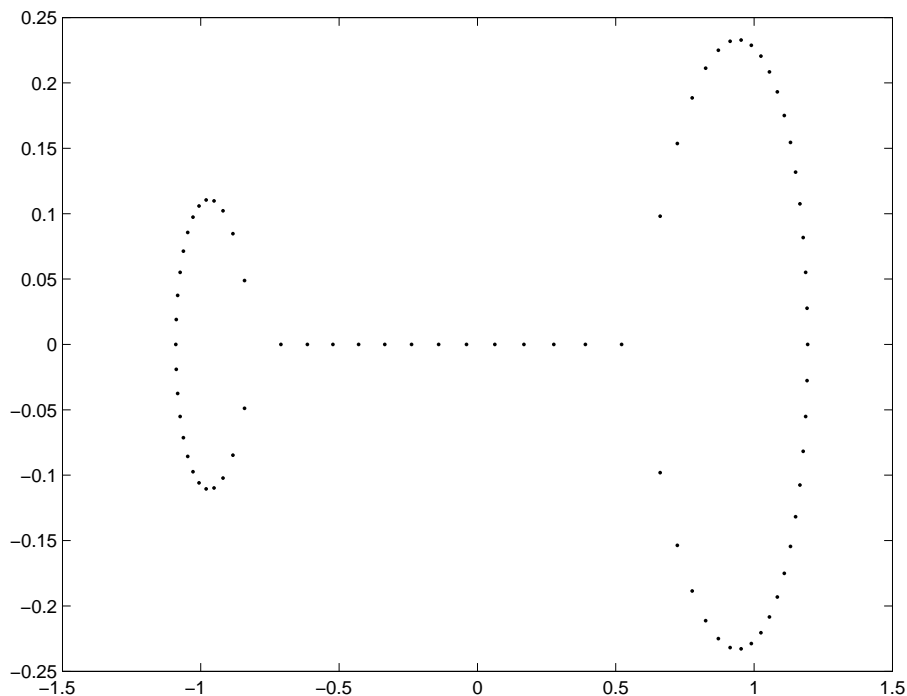}
  \caption{Zeros of $P_{75}^{(-37.4, -25.1)}$. We have
   $13=n-[-\alpha]-[-\beta]$ zeros on $(-1,1)$, $25=[-\beta]$
  zeros aligned on a curve on the left ($\Gamma_{1}$),
  and $37=[-\alpha]$ zeros aligned on a curve on the right ($\Gamma_{-1}$).}\label{fig:zerosJac}
\end{figure}
%%%%%%%%%%%%%%%%%%%%%%%%%%%%%%%%%%%%%%%%%%%%%%%%%%%%%%%%%%
%

This connection seems to go beyond the real zeros. In the example
in Fig.\ \ref{fig:zerosJac} the number of zeros on $[-1,1]$,
$\Gamma_{-1}$, and $\Gamma_1$ matches the number of orthogonality
conditions on the corresponding curves, given by Theorem
\ref{prop54}. We cannot prove this experimental fact, but we
expect to be able to establish this in an asymptotic sense by
means of the techniques from
\cite{Kuijlaars/Mclaughlin:01a}--\cite{MR1858305}.

\section*{Acknowledgements}

Arno Kuijlaars and Andrei Mart{\'\i}nez-Finkelshtein thank the organizers
of the International Workshop on Orthogonal Polynomials IWOP'02
for the invitation to speak at this workshop.

The research of A.B.J.K.\ and A.M.F.\ was partially supported by
the INTAS project 2000--272, and by the Ministry of Science and
Technology (MCYT) of Spain and the European Regional Development
Fund (ERDF) through the grant BFM2001-3878-C02-02. A.B.J.K.\ was
also supported by FWO research project G.0176.02. Additionally,
A.M.F.\ acknowledges the support of Junta de Andaluc{\'\i}a, Grupo de
Investigaci{\'o}n FQM 0229.
The research of R.O.\ was partially supported by Research Projects of
Spanish D.G.I.C.Y.T.\ and Gobierno Aut\'onomo de Canarias, under contracts
BFM2001-3411 and PI2002/136, respectively.

%\bibliographystyle{siam}
%\bibliography{../comun/jacobi,../comun/sobolev_new,../comun/general,../comun/mistrabajos,../comun/rh}

\end{document}